\documentclass[a4paper,twoside,10pt]{article}

\usepackage[top=30mm,bottom=30mm,left=30mm,right=30mm]{geometry}
\usepackage[TS1,T1]{fontenc}
\usepackage{amsmath}
\usepackage{amsfonts}
\usepackage{amssymb}
\usepackage{amsthm}
\usepackage{mathtools}
\usepackage{titling}
\usepackage{calc}
\usepackage{relsize}
\usepackage{pdflscape}
\usepackage{url}
\usepackage{hyperref}
\usepackage[USenglish=nohyphenation]{hyphsubst}
\usepackage[USenglish]{babel}
\usepackage[round]{natbib}

\newcommand{\npar}{\vspace{1.6ex plus 4pt}\par}
\newcommand{\adnpar}{\vspace{1.6ex plus 4pt}\par}

\newcommand{\Fsharp}{{\settoheight{\dimen0}{C}F\kern-.05em\resizebox{!}{\dimen0}{\raisebox{0.84\depth}{\#}}}}

\makeatletter

\g@addto@macro\normalsize{
\setlength\abovedisplayskip{0.8ex plus 3pt}
\setlength\belowdisplayskip{0ex plus 1pt}
\setlength\abovedisplayshortskip{0pt plus 1pt}
\setlength\belowdisplayshortskip{0pt plus 1pt}
}
\renewcommand\section{\@startsection{section}{1}{\z@}{-7ex \@plus -1ex \@minus -.2ex}{2ex \@plus.2ex}
{\normalfont\large\scshape\centering}}
\renewcommand\subsection{\@startsection{subsection}{1}{\z@}{-5ex \@plus -1ex \@minus -.2ex}{2ex \@plus.2ex}
{\normalfont\itshape\centering}}

\makeatother

\pretitle{\vspace*{-40.19449pt plus -3.0pt minus -5.0pt}\begin{center}\LARGE}
\posttitle{\end{center}}
\preauthor{\vspace*{0pt plus 0.66667fill}\vspace*{-20.19444pt plus -3.0pt minus -5.0pt}\vspace*{-5.384pt}\begin{center}\large}
\postauthor{\end{center}\vspace*{\fill}\vspace*{-36.75002pt plus -6.0pt minus -10.0pt}\vspace*{-3.65pt}}
\predate{}
\postdate{}

\def\uprighte{\mathrm{e}}
\def\uprighti{\mathrm{i}}
\def\uprightpi{\mathrm{\pi}}
\def\uprightpartial{\reflectbox{$6$}}

\def\uprightpartialspace{\uprightpartial\:\!}
\def\uprightpartialoperator#1{\dfrac{\uprightpartial}{\uprightpartialspace{#1}}}
\def\nuprightpartialoperator#1#2{\dfrac{\uprightpartial^{#1}}{\uprightpartialspace{#2}}}
\def\htext#1{\text{{#1}\vphantom{h}}}
\def\hvar#1{{#1}\vphantom{h}}

\def\matrixwidth#1{\mathmakebox[\widthof{60}][c]{{#1}}}

\def\oneminusoneone{{\def\arraystretch{0.5}\biggl\{\begin{matrix}1\\-1\\1\end{matrix}\biggr\}}}
\def\oneminusoneoneone{{\def\arraystretch{0.5}\Biggl\{\begin{matrix}1\\-1\\1\\1\end{matrix}\Biggr\}}}
\def\minusoneoneminusoneone{{\def\arraystretch{0.5}\Biggl\{\begin{matrix}-1\\1\\-1\\1\end{matrix}\Biggr\}}}
\def\minusoneoneoneone{{\def\arraystretch{0.5}\Biggl\{\begin{matrix}-1\\1\\1\\1\end{matrix}\Biggr\}}}
\def\iioneone{{\def\arraystretch{0.5}\Biggl\{\begin{matrix}\uprighti\\\uprighti\\1\\1\end{matrix}\Biggr\}}}
\def\ioneoneone{{\def\arraystretch{0.5}\Biggl\{\begin{matrix}{\uprighti}\\1\\1\\1\end{matrix}\Biggr\}}}
\def\smalliioneone{{\def\arraystretch{0.5}\biggl\{\begin{smallmatrix}\uprighti\\\uprighti\\1\\1\end{smallmatrix}\biggr\}}}
\def\twoitwotwotwo{{\def\arraystretch{0.5}\Biggl\{\begin{matrix}2\, \uprighti\\2\\2\\2\end{matrix}\Biggr\}}}
\def\sincossinhcosh{{\def\arraystretch{0.5}\Biggl\{\begin{matrix}\htext{sin}\\\htext{cos}\\\text{sinh}\\\text{cosh}\end{matrix}\Biggr\}}}
\def\cossincoshsinh{{\def\arraystretch{0.5}\Biggl\{\begin{matrix}\htext{cos}\\\htext{sin}\\\text{cosh}\\\text{sinh}\end{matrix}\Biggr\}}}
\def\smallcossincoshsinh{{\def\arraystretch{0.5}\biggl\{\begin{smallmatrix}\htext{cos}\\\htext{sin}\\\text{cosh}\\\text{sinh}\end{smallmatrix}\biggr\}}}
\def\fgh{{\def\arraystretch{0.5}\biggl\{\begin{matrix}\hvar{f}\\\hvar{g}\\\hvar{h}\end{matrix}\biggr\}}}
\def\fgfh{{\def\arraystretch{0.5}\Biggl\{\begin{matrix}\hvar{f}\\\hvar{g}\\\hvar{f}\\\hvar{h}\end{matrix}\Biggr\}}}
\def\sctd{{\def\arraystretch{0.5}\Biggl\{\begin{matrix}\hvar{s}\\\hvar{c}\\\hvar{t}\\\hvar{d}\end{matrix}\Biggr\}}}
\def\complexexpression{{\def\arraystretch{0.5}\left\{\begin{matrix}\uprighte^{\uprighti x} - \uprighte^{-\uprighti x}\\
\uprighte^{\uprighti x} + \uprighte^{-\uprighti x}\\
\uprighte^x - \uprighte^{-x}\\
\uprighte^x + \uprighte^{-x}\end{matrix}\right\}}}
\def\largesum{\mathlarger{\mathlarger{\mathlarger{\sum}}}}
\def\lftse{\hspace{-0.3em}}
\def\lftfp{\hspace{-0.263em}}
\def\powersmatrix{\left[\setlength{\arraycolsep}{0.25em}\def\arraystretch{0.5}\begin{array}{cccccc}
\matrixwidth{1} & \matrixwidth{} & \matrixwidth{}
& \matrixwidth{} & \matrixwidth{} & \matrixwidth{} \\
  &  1 &    &    &    &    \\
  &    &  2 &    &    &    \\
  &    &    &  4 &    &    \\
  &    &    &    &  8 &    \\
  &    &    &    &    & 16 \\
\multicolumn{6}{r}{\vdots}
\end{array}\right]}
\def\heightened#1{\mathmakebox[\widthof{$\powersmatrix$}][c]{\textrm{\raisebox{29.42pt}{#1}}}}
\def\derivativesleft{\nuprightpartialoperator{k}{x^k}\left[\sincossinhcosh^{\lftse n}(x)\right]}
\def\derivativesright{\dfrac{\minusoneoneminusoneone^{\lftse n}}{2^n}
\cdot \largesum_{r = 0}^n\,\minusoneoneminusoneone^{\lftse r}
\cdot \dbinom{n}{r}
\cdot \ioneoneone^{\lftse -n}
\cdot \iioneone^k
\cdot (2 r - n)^k
\cdot \uprighte^{\smalliioneone \cdot (2 r - n) \cdot x}.}
\def\exceptzeroes{{\left\{\def\arraystretch{0.5}\left\{\begin{matrix}m\\(m + \tfrac{1}{2})\\\uprighti\,m\\\uprighti\,(m + \tfrac{1}{2})\end{matrix}\right\}\uprightpi\,|\,n < k \wedge m \in \mathbb{Z}\right\}}}
\def\resultbox#1{\mathmakebox[\widthof{$\nuprightpartialoperator{0}{x^0} \bigl[\cosh^n(x)\bigr]$}][r]{#1}}
\def\resultleft#1#2{\resultbox{\nuprightpartialoperator{#2}{x^{#2}} \bigl[{#1}^n(x)\bigr]}}
\def\pnl{\phantom{\Bigr]}}
\def\inductionhypothesis{\begin{equation*}
\nuprightpartialoperator{k}{x^k}\left[\sincossinhcosh^{\lftse n}(x)\right]
= \sincossinhcosh^{\lftse n - k}(x) \cdot \fgfh_{\lftse k}(n)\left(\cossincoshsinh(x)\right)\lftfp.
\end{equation*}}

\def\arraystretch{0.5}

\defcitealias{oeis}{Adamson~(2007)}
\defcitealias{knuth}{Knuth~(1992)}
\defcitealias{pascal}{Pascal's triangle~(1655)}
\defcitealias{qi}{Qi~(2015)}

\title{On the derivatives of the powers\\*
of trigonometric and hyperbolic sine and cosine}
\author{Stijn ``\kern -0.16em Adhemar\kern 0.04em '' Vandamme}
\date{}

\begin{document}
\pagestyle{empty}
\maketitle
\thispagestyle{empty}
\vspace{1em}
\begin{abstract}
\vspace{1.32494pt}
\noindent
{\sloppy This work contains different expressions for the $k$th derivative of the $n$th power of the \linebreak[4]
trigonometric and hyperbolic sine and cosine.}
\vspace{6pt plus 3.0pt minus 6.0pt}
\par
\noindent
{\fussy The first set of expressions follow from the complex definitions of the trigonometric and hyperbolic sine and cosine, and the binomial theorem.}
\vspace{6pt plus 3.0pt minus 6.0pt}
\par
\noindent
The other expressions are polynomial-based. They are perhaps less obvious, and use only polynomials in $\sin x$ and $\cos x$, or in $\sinh x$ and $\cosh x$.
No sines or cosines of arguments other than $x$ appear in these polynomial-based expressions.
The final expressions are dependent only on $\sin x$, $\cos x$, $\sinh x$, or $\cosh x$ respectively when $k$ is even;
and they only have a single additional factor $\cos x$, $\sin x$, $\cosh x$, or $\sinh x$ respectively when $k$ is odd.
\end{abstract}
\vspace*{\fill}
\vspace*{-6pt plus -3.0pt minus -5.0pt}
\vspace*{-5.16667pt}
\vspace*{-30.1388pt plus -4.30554pt minus -0.86108pt}
\section*{The derivatives of the powers of the natural exponential function}
\label{exponential}
This works starts with expressions for the derivatives of the powers of the natural exponential function.
\npar
For $k \in \mathbb{N}$, and $n \in \mathbb{Z}$, the $k$th derivative of the function applying the $n$th power to the natural exponential is:
\begin{equation*}
\nuprightpartialoperator{k}{x^k} \bigl[\exp^n x\bigr] = n^k \exp^n x = n^k \exp {(n x)} = n^k \uprighte^{n x}.
\end{equation*}
\adnpar
The natural exponential function, its powers, and their derivatives are all both $\mathbb{C} \rightarrow \mathbb{C}$ and $\mathbb{R} \rightarrow \mathbb{R}$ functions:
they are defined for the whole domain of the complex numbers, returning complex values; yet the image of all real arguments are themselves real values.
\subsection*{The definedness of exponentation}
\label{exponentation}
For the purposed of this work, $0^0$ is considered defined as $0^0 = 1$. For a short discussion and historical overview of this controversy, see \citetalias{knuth}.
So the following identity holds over the entire domain $\mathbb{C}$ (including $\mathbb{R}$, including $\{0\}$):
\begin{equation*}
\exp^0 x = \nuprightpartialoperator{0}{x^0} \bigl[\exp^0 x\bigr] = 0^0\exp^0 x = 0^0\exp 0 = 1.
\end{equation*}
\adnpar
For bases $z \in \mathbb{C}$, the power $z^n$ is defined for all exponents $n \in \mathbb{Z}$.
For bases $z \in \mathbb{R}_{\geqslant 0}$, the power $z^n$ is also defined for all exponents $n \in \mathbb{R}$.
As the power rule also applies for those non-integer exponents,
the polynomial-based identities obtained in this work also apply for real powers,
in the intervals for $x$ in which $\sin x$, $\cos x$, $\sinh x$, or $\cosh x$ is real and nonnegative.
Note that $\sinh x \in \mathbb{R}_{\geqslant 0}$ when $x \in \mathbb{R}_{\geqslant 0}$, and that $\cosh x \in \mathbb{R}_{\geqslant 0}$ when $x \in \mathbb{R}$.
However, the 4 functions are not all
\linebreak[2]nonnegative-real-valued (or even real-valued) for every $z \in \mathbb{C}$.
In light of this note, only the (strictly too strict) condition $n \in \mathbb{Z}$ will be considered and mentioned in this work.
\section*{Using the complex definitions and the binomial theorem}
\label{binomial}
The first set of expressions for the derivatives of the powers of trigonometric and hyperbolic sine and cosine
follow from the complex definitions of the trigonometric and hyperbolic sine and cosine, which are:
\begin{equation*}
\sin x = \frac{\uprighte^{\uprighti x} - \uprighte^{-\uprighti x}}{2\, \uprighti},
\quad\quad
\cos x = \frac{\uprighte^{\uprighti x} + \uprighte^{-\uprighti x}}{2},
\quad\quad
\sinh x = \frac{\uprighte^x - \uprighte^{-x}}{2},
\quad\quad
\cosh x = \frac{\uprighte^x + \uprighte^{-x}}{2}.
\end{equation*}
\adnpar
These functions are all both $\mathbb{C} \rightarrow \mathbb{C}$ and $\mathbb{R} \rightarrow \mathbb{R}$ functions.
\npar
It follows from these complex definitions that powers of the trigonometric and hyperbolic sine and cosine are:
\begin{equation*}
\setlength{\arraycolsep}{1.4pt}
\begin{array}{r c l}
\mathmakebox[\widthof{$\derivativesleft$}][r]{\sincossinhcosh^{\lftse n}(x)}
& = &
\mathmakebox[\widthof{$\derivativesright$}][l]{\dfrac{1}{{\twoitwotwotwo}^{\lftse n}} \cdot \complexexpression^{\lftse n}.}
\end{array}
\end{equation*}
\adnpar
For $n \in \mathbb{N}$, the binomial theorem can be applied:
\begin{equation*}
\setlength{\arraycolsep}{1.4pt}
\begin{array}{r c l}
\mathmakebox[\widthof{$\derivativesleft$}][r]{\sincossinhcosh^{\lftse n}(x)}
& = &
\mathmakebox[\widthof{$\derivativesright$}][l]{\dfrac{1}{{\twoitwotwotwo}^{\lftse n}}
\cdot \largesum_{r = 0}^n\,\minusoneoneminusoneone^{\lftse n - r}
\cdot \dbinom{n}{r}
\cdot \uprighte^{\smalliioneone \cdot r \cdot x}
\cdot \uprighte^{-\smalliioneone \cdot (n - r) \cdot x}}\\
\\
& = &
\dfrac{\minusoneoneminusoneone^{\lftse n}}{{\twoitwotwotwo}^{\lftse n}}
\cdot \largesum_{r = 0}^n\,\minusoneoneminusoneone^{\lftse r}
\cdot \dbinom{n}{r}
\cdot \uprighte^{\smalliioneone \cdot (2 r - n) \cdot x}\\
\\
& = &
\dfrac{\minusoneoneminusoneone^{\lftse n}}{2^n}
\cdot \largesum_{r = 0}^n\,\minusoneoneminusoneone^{\lftse r}
\cdot \dbinom{n}{r}
\cdot \ioneoneone^{\lftse -n}
\cdot \uprighte^{\smalliioneone \cdot (2 r - n) \cdot x}.
\end{array}
\end{equation*}
\adnpar
And so the derivatives are:
\begin{equation*}
\setlength{\arraycolsep}{1.4pt}
\begin{array}{r c l}
\derivativesleft
& = &
\derivativesright
\end{array}
\end{equation*}
\adnpar
For trigonometric sine:
\begin{equation*}
\setlength{\arraycolsep}{1.4pt}
\begin{array}{r r c l l}
&\mathmakebox[\widthof{$\derivativesleft$}][r]{\nuprightpartialoperator{k}{x^k}\bigl[\sin^n x\bigr]}
& = &
\mathmakebox[\widthof{$\derivativesright$}][l]{\dfrac{(-1)^n}{2^n}
\cdot \mathlarger{\sum}_{r = 0}^n\,(-1)^r
\cdot \dbinom{n}{r}
\cdot (2 r - n)^k
\cdot \uprighti^{k - n}
\cdot \uprighte^{\uprighti (2 r - n) x}}
&\\
\\
& & = &
\dfrac{(-1)^n}{2^n}
\cdot \mathlarger{\sum}_{r = 0}^n\,(-1)^r
\cdot \dbinom{n}{r}
\cdot (2 r - n)^k
\cdot \uprighte^{\uprighti\,\left((k - n)\frac{\pi}{2} + (2 r - n) x\right)} &\\
\\
\multicolumn{5}{r}{
=
\dfrac{(-1)^n}{2^n}
\cdot \mathlarger{\sum}_{r = 0}^n\,(-1)^r
\cdot \dbinom{n}{r}
\cdot (2 r - n)^k
\cdot \Bigl[\cos\bigl(\frac{(k - n) \pi}{2} + (2 r - n) x\bigr)
+ \uprighti\sin\bigl(\frac{(k - n) \pi}{2} + (2 r - n) x\bigr)\Bigr].}
\end{array}
\end{equation*}
\adnpar
For trigonometric cosine:
\begin{equation*}
\setlength{\arraycolsep}{1.4pt}
\begin{array}{r r c l l}
&\mathmakebox[\widthof{$\derivativesleft$}][r]{\nuprightpartialoperator{k}{x^k}\bigl[\cos^n x\bigr]}
& = &
\mathmakebox[\widthof{$\derivativesright$}][l]{\dfrac{1}{2^n}
\cdot \mathlarger{\sum}_{r = 0}^n\,\dbinom{n}{r}
\cdot (2 r - n)^k
\cdot \uprighti^k
\cdot \uprighte^{\uprighti (2 r - n) x}}
&\\
\\
& & = &
\dfrac{1}{2^n}
\cdot \mathlarger{\sum}_{r = 0}^n\,\dbinom{n}{r}
\cdot (2 r - n)^k
\cdot \uprighte^{\uprighti\,\left(\frac{k \pi}{2} + (2 r - n) x\right)} &\\
\\
& & = &
\dfrac{1}{2^n}
\cdot \mathlarger{\sum}_{r = 0}^n\,\dbinom{n}{r}
\cdot (2 r - n)^k
\cdot \Bigl[\cos\bigl(\frac{k \pi}{2} + (2 r - n) x\bigr)
+ \uprighti\sin\bigl(\frac{k \pi}{2} + (2 r - n) x\bigr)\Bigr]. &
\end{array}
\end{equation*}
\adnpar
In the above results for trigonometric sine and cosine, the imaginary part evaluates to 0 for real arguments $x$.
\npar
The above results for trigonometric sine and cosine are not new, as they are already published
\linebreak[1] by \citetalias{qi}.
\npar
For hyperbolic sine:
\begin{equation*}
\setlength{\arraycolsep}{1.4pt}
\begin{array}{r r c l l}
&\mathmakebox[\widthof{$\derivativesleft$}][r]{\nuprightpartialoperator{k}{x^k}\bigl[\sinh^n x\bigr]}
& = &
\mathmakebox[\widthof{$\derivativesright$}][l]{\dfrac{(-1)^n}{2^n}
\cdot \mathlarger{\sum}_{r = 0}^n\,(-1)^r
\cdot \dbinom{n}{r}
\cdot (2 r - n)^k
\cdot \uprighte^{(2 r - n) x}}\\
\\
& & = &
\dfrac{(-1)^n}{2^n}
\cdot \mathlarger{\sum}_{r = 0}^n\,(-1)^r
\cdot \dbinom{n}{r}
\cdot (2 r - n)^k
\cdot \Bigl[\cosh\bigl((2 r - n) x\bigr) + \sinh\bigl((2 r - n) x\bigr)\Bigr]. &
\end{array}
\end{equation*}
\adnpar
For hyperbolic cosine:
\begin{equation*}
\setlength{\arraycolsep}{1.4pt}
\begin{array}{r r c l l}
&\mathmakebox[\widthof{$\derivativesleft$}][r]{\nuprightpartialoperator{k}{x^k}\bigl[\cosh^n x\bigr]}
& = &
\mathmakebox[\widthof{$\derivativesright$}][l]{\dfrac{1}{2^n}
\cdot \mathlarger{\sum}_{r = 0}^n\,\dbinom{n}{r}
\cdot (2 r - n)^k
\cdot \uprighte^{(2 r - n) x}}
&\\
\\
& & = &
\dfrac{1}{2^n}
\cdot \mathlarger{\sum}_{r = 0}^n\,\dbinom{n}{r}
\cdot (2 r - n)^k
\cdot \Bigl[\cosh\bigl((2 r - n) x\bigr) + \sinh\bigl((2 r - n) x\bigr)\Bigr]. &
\end{array}
\end{equation*}
\section*{Using polynomials, intermediate step}
\label{intermediate}
\subsection*{Definitions of the intermediate polynomials}
\label{intermediatepolynomials}
For $k \in \mathbb{N}$, $n \in \mathbb{Z}$, and $u \in \mathbb{C}$,
let $f_k$, $g_k$, and $h_k$
be sequences of curried polynomial functions, whose images
$f_k(n)(u)$, $g_k(n)(u)$, and $h_k(n)(u)$
are polynomials in $u$, with coefficients that are,
\linebreak[1]in turn, polynomials in $n$;
and let those sequences be recursively defined (over $k$) as follows:
\begin{equation*}
\setlength{\arraycolsep}{1.4pt}
\begin{array}{r c l}
\fgh_{\lftse 0}(n)(u) & = & 1,\\
\\
\fgh_{\lftse k + 1}(n)(u) & = & \oneminusoneone
\cdot \Biggl((n-k)
\cdot u
\cdot \fgh_{\lftse k}(n)(u)
+
\left\{{\def\arraystretch{0.5}\begin{matrix}u^2-1\\u^2-1\\u^2+1\end{matrix}}\right\}
\cdot \uprightpartialoperator{u} \left[\fgh_{\lftse k}(n)(u)\right]\Biggr).
\end{array}
\end{equation*}
\subsection*{Properties}
\label{intermediateproperties}
The rank of the polynomials $f_k(n)(u)$, $g_k(n)(u)$, and $h_k(n)(u)$ is $k$.
\npar
The coefficient of the highest power with non-zero coefficient of $f_k(n)$ and of $h_k(n)$ is $1$.
The coefficient of the highest power with non-zero coefficient of $g_k(n)$ is $1$ when $k$ is even, or $-1$ when $k$ is odd.
\npar
The polynomials $f_k(n)(u)$, $g_k(n)(u)$, and $h_k(n)(u)$ contain only non-zero coefficients for even powers of $u$ when $k$ is even,
or only non-zero coefficients for odd powers of $u$ when $k$ is odd.
\npar
For $k \in \mathbb{N}$, $r \in \mathbb{N}$, when $r \leqslant k$, and when $k$ and $r$ have the same parity, the coefficients for $u^r$ are polynomials in $n$ with rank $\dfrac{k + r}{2}$.
\npar
$\forall\,k \in \mathbb{N}$, $n \in \mathbb{Z}$, and $u \in \mathbb{C}$: $g_{k}(n)(u) = (-1)^k f_{k}(n)(u).$
\npar
The coefficients of the second-highest power with non-zero coefficient of $f_{k}(n)$ and of $g_{k}(n)$ are polynomials in $n$,
of which the signs of the coefficients are alternating,
and the absolute values of the coefficients appear in the lower-triangular matrix obtained through the multiplication of
the following \nolinebreak[4]3 \nolinebreak[4]lower-triangular matrices,
followed by the removal of the left column of the
\linebreak[4]lower-triangular product matrix, which is typeset in bold below:
\begin{equation*}
\begin{array}{c}
\left[\setlength{\arraycolsep}{0.25em}\def\arraystretch{0.5}\begin{array}{cccccc}
\matrixwidth{1} & \matrixwidth{} & \matrixwidth{}
& \matrixwidth{} & \matrixwidth{} & \matrixwidth{} \\
1 &  1 &    &    &    &    \\
1 &  2 &  1 &    &    &    \\
1 &  3 &  3 &  1 &    &    \\
1 &  4 &  6 &  4 &  1 &    \\
1 &  5 & 10 & 10 &  5 &  1 \\
\multicolumn{6}{c}{\vdots}
\end{array}\right]
\times
\left[\setlength{\arraycolsep}{0.25em}\def\arraystretch{0.5}\begin{array}{cccccc}
\matrixwidth{1} & \matrixwidth{} & \matrixwidth{}
& \matrixwidth{} & \matrixwidth{} & \matrixwidth{} \\
1 &  1 &    &    &    &    \\
  &  1 &  1 &    &    &    \\
  &    &  1 &  1 &    &    \\
  &    &    &  1 &  1 &    \\
  &    &    &    &  1 &  1 \\
\multicolumn{6}{r}{\vdots}
\end{array}\right]
\times
\powersmatrix\\
\\
\heightened{\citetalias{pascal}}
\mathmakebox[\widthof{=}][c]{=}
\mathmakebox[\widthof{\;}][c]{}
\left[\setlength{\arraycolsep}{0.25em}\def\arraystretch{0.5}\begin{array}{cccccc}
\matrixwidth{\textbf{1}} & \matrixwidth{} & \matrixwidth{}
& \matrixwidth{} & \matrixwidth{} & \matrixwidth{} \\
\textbf{2} &  1 &    &    &    &    \\
\textbf{3} &  3 &  2 &    &    &    \\
\textbf{4} &  6 &  8 &  2 &    &    \\
\textbf{5} & 10 & 20 & 20 &  8 &    \\
\textbf{6} & 15 & 40 & 60 & 48 & 16 \\
\multicolumn{6}{c}{\vdots}
\end{array}\right]
\mathmakebox[\widthof{=}][l]{.}
\mathmakebox[\widthof{\;}][c]{}
\heightened{Powers of 2}
\end{array}
\end{equation*}
\adnpar
The sequence of numbers in this lower-triangular product matrix, including its left column, has been contributed by \citetalias{oeis}
to the On-Line Encyclopedia of Integer Sequences (\textsc{oeis}), where it is registered as sequence A133341.
\subsection*{The derivatives of the powers of trigonometric and hyperbolic sine and cosine}
\label{intermediatederivations}
With $f_k(n)(u)$, $g_k(n)(u)$, and $h_k(n)(u)$ defined above, the following identities will be proved:
\inductionhypothesis
\subsection*{Proof}
\label{proof}
The proof is by induction on $k$.
\npar
The base case, for $k = 0$, is trivial:
\begin{equation*}
\nuprightpartialoperator{0}{x^0}\left[\sincossinhcosh^{\lftse n}(x)\right]
= \sincossinhcosh^{\lftse n}(x)
= \sincossinhcosh^{\lftse n - 0}(x) \cdot 1
= \sincossinhcosh^{\lftse n - 0}(x) \cdot \fgfh_{\lftse 0}(n)\left(\cossincoshsinh(x)\right)\lftfp.
\end{equation*}
\adnpar
For the inductive case, assume the induction hypothesis
\inductionhypothesis
\adnpar
With $f_k(n)(u)$, $g_k(n)(u)$, and $h_k(n)(u)$ defined above, under this induction hypothesis,
and using the product rule and the chain rule, the $(k + 1)$th derivatives can intermediately be derived
\linebreak[4]as follows:
\begin{landscape}
\abovedisplayskip=0pt
\belowdisplayskip=0pt
\abovedisplayshortskip=0pt
\belowdisplayshortskip=0pt
\vspace*{-0.823cm}{
\[
\def\arraystretch{0.8}\begin{array}{l}
\nuprightpartialoperator{k + 1}{x^{k + 1}}\left[\sincossinhcosh^{\lftse n}(x)\right]
\\
\setlength{\arraycolsep}{1.4pt}\def\arraystretch{4.1}\begin{array}{crcl}
= &
\multicolumn{3}{l}{
\uprightpartialoperator{x}\left[\sincossinhcosh^{\lftse n - k}(x)
\cdot \fgfh_{\lftse k}(n)\left(\cossincoshsinh(x)\right)\right]}
\\
= &
\uprightpartialoperator{x}\left[\sincossinhcosh^{\lftse n - k}(x)\right]
\cdot \fgfh_{\lftse k}(n)\left(\cossincoshsinh(x)\right)
& + &
\sincossinhcosh^{\lftse n - k}(x)
\cdot \uprightpartialoperator{x}\left[\fgfh_{\lftse k}(n)\left(\cossincoshsinh(x)
\right)\right]
\\
= &
(n - k)
\cdot \sincossinhcosh^{\lftse n - k - 1}(x)
\cdot \uprightpartialoperator{x}\left[\sincossinhcosh(x)\right]
\cdot \fgfh_{\lftse k}(n)\left(\cossincoshsinh(x)\right)
& + &
\sincossinhcosh^{\lftse n - k - 1}(x)
\cdot \sincossinhcosh(x)
\cdot \uprightpartialoperator{u}\left[\fgfh_{\lftse k}(n)(u)\right]_{u = \smallcossincoshsinh(x)}
\\
\multicolumn{4}{r}{
\cdot \uprightpartialoperator{x}\left[\left(\cossincoshsinh(x)\right)
\right]}
\\
= &
\sincossinhcosh^{\lftse n - k - 1}(x)
\cdot (n - k)
\cdot \oneminusoneoneone
\cdot \cossincoshsinh(x)
\cdot \fgfh_{\lftse k}(n)\left(\cossincoshsinh(x)\right)
& + &
\sincossinhcosh^{\lftse n - k - 1}(x)
\cdot \minusoneoneoneone
\cdot \sincossinhcosh^{\lftse 2}(x)
\cdot \uprightpartialoperator{u}\left[\fgfh_{\lftse k}(n)(u)\right]_{u = \smallcossincoshsinh(x)}
\\
= &
\sincossinhcosh^{\lftse n - k - 1}(x)
\cdot \oneminusoneoneone
\cdot \left(\vphantom{\uprightpartialoperator{x}\left[\fgfh_{\lftse k}(n)(u)\right]_{u = \smallcossincoshsinh(x)}}\right.\hspace{-0.35em}(n - k)
\cdot \cossincoshsinh(x)
\cdot \fgfh_{\lftse k}(n)\left(\cossincoshsinh(x)\right)
& + &
\left.{\hspace{-0.12em}\left\{\setlength{\arraycolsep}{1.4pt}\def\arraystretch{0.5}\begin{matrix*}[r]\cos^2 x&\mathop{-} 1\\
\sin^2 x&\mathop{-} 1\\
\cosh^2 x&\mathop{-} 1\\
\sinh^2 x&\mathop{+} 1\end{matrix*}\right\}}
\cdot \uprightpartialoperator{u}\left[\fgfh_{\lftse k}(n)(u)\right]_{u = \smallcossincoshsinh(x)}\right)
\\
= &
\multicolumn{2}{l}{
\sincossinhcosh^{\lftse n - (k + 1)}(x)
\cdot \fgfh_{\lftse k + 1}(n)\left(\cossincoshsinh(x)\right)\lftfp.} &
\multicolumn{1}{r}{\raisebox{-13.082pt}{$\square\mathmakebox[0.75pt][l]{}$}}
\end{array}
\end{array}
\]
}
\end{landscape}
\section*{Using polynomials, final step}
\label{final}
\subsection*{Definitions of the final polynomials}
\label{finalpolynomials}
For $k \in \mathbb{N}$, $n \in \mathbb{Z}$, and $u \in \mathbb{C}$,
let $s_k$, $c_k$, $t_k$, and $d_k$
be sequences of curried polynomial functions,
whose images $f_k(n)(u)$, $g_k(n)(u)$, and $h_k(n)(u)$
are polynomials in $u$, with coefficients that are,
\linebreak[1]in turn, polynomials in $n$;
and let those sequences be obtained as follows:
\begin{equation*}
\setlength{\arraycolsep}{1.4pt}
\begin{array}{rcl}
\sctd_{\lftse k}(n)(v)
& = &
\textrm{substitute (for all $m \in \mathbb{N}$) $u^{2m}$ with }{\def\arraystretch{0.5}\left\{\begin{matrix}1 - v^{2m}\\
1 - v^{2m}\\
v^{2m} + 1\\
v^{2m} - 1\end{matrix}\right\}}\\
\\
& & \textrm{in }\fgfh_{\lftse k}(n)(u)\textrm{ when $k$ is even, or in }\dfrac{1}{u} \fgfh_{\lftse k}(n)(u)\textrm{ when $k$ is odd.}
\end{array}
\end{equation*}
\subsection*{Properties}
\label{finalproperties}
The rank of the polynomials $s_k(n)(v)$, $c_k(n)(v)$, $t_k(n)(v)$, and $d_k(n)(v)$ is always even:
it $k$ when $k$ is even, or $k - 1$ when $k$ is odd.
\npar
The coefficient of the highest power with non-zero coefficient of $t_k(n)(v)$ and of $d_k(n)(v)$ is $1$.
The coefficient of the highest power with non-zero coefficient of $s_k(n)(v)$ is $1$ when $k \bmod 4 \in \{0,\,1\}$, or $-1$ when $k \bmod 4 \in \{2,\,3\}$.
The coefficient of the highest power with non-zero coefficient
\linebreak[2]of $c_k(n)(v)$ is $1$ when $k \bmod 4 \in \{0,\,3\}$, or $-1$ when $k \bmod 4 \in \{1,\,2\}$.
\npar
The polynomials $s_k(n)(v)$, $c_k(n)(v)$, $t_k(n)(v)$ and $d_k(n)(v)$ contain only non-zero coefficients
\linebreak[4]for even powers of $v$.
\npar
$\forall\,k \in \mathbb{N}$, $n \in \mathbb{Z}$, and $v \in \mathbb{C}$: $c_{k}(n)(v) = (-1)^k s_{k}(n)(v).$
\subsection*{The derivatives of the powers of trigonometric and hyperbolic sine and cosine}
With $s_k(n)(u)$, $c_k(n)(u)$, $t_k(n)(u)$ and $d_k(n)(u)$ defined above, and using the identities $\forall\,x \in \mathbb{C}: \sin^2(x) + \cos^2(x) = 1$ and $\cosh^2(x) - \sinh^2(x) = 1$, the $k$th derivatives can finally be derived
\linebreak[1]as follows:
\begin{equation*}
\nuprightpartialoperator{k}{x^k}\left[\sincossinhcosh^{\lftse n}(x)\right] =
\left\{\setlength{\arraycolsep}{1.4pt}
\begin{array}{ll}
\sincossinhcosh^{\lftse n - k}(x) \cdot \sctd_{\lftse k}(n)\left(\sincossinhcosh(x)\right)
& \textrm{ when $k$ is even,}\\
\\
\sincossinhcosh^{\lftse n - k}(x) \cdot \cossincoshsinh(x) \cdot \sctd_{\lftse k}(n)\left(\sincossinhcosh(x)\right)
& \textrm{ when $k$ is odd.}
\end{array}\right.
\end{equation*}
\section*{Using polynomials, (partially) applied results}
\label{results}
Both the intermediate and the final expressions are listed here.
The expression in brackets is the intermediate or final polynomial, except for when the coefficient of the highest power with non-zero coefficient of $g_k(n)$, of $s_k(n)$, or of $c_k(n)$ is $-1$, in which case $-1$ is factored out.\phantom{\citet{hoffman}}
\npar
For $k \in \mathbb{N}$, $n \in \mathbb{Z}$, and $x \in \mathbb{C} \setminus \exceptzeroes$:
\subsection*{The derivatives of the powers of trigonometric sine}
\label{sin}
\vspace*{-0.8ex}
$\resultleft{\sin}{0} = \sin^{n-0}(x) \bigl[1\bigr]\pnl$\\
$\resultbox{} = \sin^{n-0}(x) \bigl[1\bigr] = \sin^n(x),$\\
$\resultleft{\sin}{1} = \sin^{n-1}(x) \bigl[n \cos(x)\bigr]\pnl$\\
$\resultbox{} = \sin^{n-1}(x) \cos(x) \bigl[n\bigr],\pnl$\\
$\resultleft{\sin}{2} = \sin^{n-2}(x) \bigl[n^2 \cos^2(x) - n\bigr]\pnl$\\
$\resultbox{} = -\sin^{n-2}(x) \bigl[n^2 \sin^2(x) + (-n^2 + n)\bigr],\pnl$\\
$\resultleft{\sin}{3} = \sin^{n-3}(x) \bigl[n^3 \cos^3(x) + (-3 n^2 + 2 n) \cos(x)\bigr]\pnl$\\
$\resultbox{} = -\sin^{n-3}(x) \cos(x) \bigl[n^3 \sin^2(x) + (-n^3 + 3 n^2 - 2 n)\bigr],\pnl$\\
$\resultleft{\sin}{4} = \sin^{n-4}(x) \bigl[n^4 \cos^4(x) + (-6 n^3 + 8 n^2 - 4 n) \cos^2(x) + (3 n^2 - 2 n)\bigr]\pnl$\\
$\resultbox{} = \sin^{n-4}(x) \bigl[n^4 \sin^4(x) + (-2 n^4 + 6 n^3 - 8 n^2 + 4 n) \sin^2(x) + (n^4 - 6 n^3 + 11 n^2 - 6 n)\bigr],\pnl$\\
$\resultleft{\sin}{5} = \sin^{n-5}(x) \bigl[n^5 \cos^5(x) + (-10 n^4 + 20 n^3 - 20 n^2 + 8 n) \cos^3(x)\pnl$\\
$\resultbox{}\phantom{ = \sin^{n-5}(x)}\;+ (15 n^3 - 30 n^2 + 16 n) \cos(x)\bigr]\pnl$\\
$\resultbox{} = \sin^{n-5}(x) \cos(x) \bigl[n^5 \sin^4(x) + (-2 n^5 + 10 n^4 - 20 n^3 + 20 n^2 - 8 n) \sin^2(x)\pnl$\\
$\resultbox{}\phantom{ = \sin^{n-5}(x) \cos(x)}\;+ (n^5 - 10 n^4 + 35 n^3 - 50 n^2 + 24 n)\bigr],\pnl$\\
$\resultleft{\sin}{6} = \sin^{n-6}(x) \bigl[n^6 \cos^6(x) + (-15 n^5 + 40 n^4 - 60 n^3 + 48 n^2 - 16 n) \cos^4(x)\pnl$\\
$\resultbox{}\phantom{ = \sin^{n-6}(x)}\;+ (45 n^4 - 150 n^3 + 196 n^2 - 88 n) \cos^2(x) + (-15 n^3 + 30 n^2 - 16 n)\bigr]\pnl$\\
$\resultbox{} = -\sin^{n-6}(x) \bigl[n^6 \sin^6(x) + (-3 n^6 + 15 n^5 - 40 n^4 + 60 n^3 - 48 n^2 + 16 n) \sin^4(x)\pnl$\\
$\resultbox{}\phantom{ = -\sin^{n-6}(x)}\;+ (3 n^6 - 30 n^5 + 125 n^4 - 270 n^3 + 292 n^2 - 120 n) \sin^2(x)\pnl$\\
$\resultbox{}\phantom{ = -\sin^{n-6}(x)}\;+ (-n^6 + 15 n^5 - 85 n^4 + 225 n^3 - 274 n^2 + 120 n)\bigr],~\ldots\pnl$
\vspace*{-2.4ex}
\subsection*{The derivatives of the powers of  trigonometric cosine}
\label{cos}
\vspace*{-0.8ex}
$\resultleft{\cos}{0} = \cos^{n-0}(x) \bigl[1\bigr]\pnl$\\
$\resultbox{} = \cos^{n-0}(x) \bigl[1\bigr] = \cos^n(x),\pnl$\\
$\resultleft{\cos}{1} = -\cos^{n-1}(x) \bigl[n \sin(x)\bigr]\pnl$\\
$\resultbox{} = -\cos^{n-1}(x) \sin(x) \bigl[n\bigr],\pnl$\\
$\resultleft{\cos}{2} = \cos^{n-2}(x) \bigl[n^2 \sin^2(x) - n\bigr]\pnl$\\
$\resultbox{} = -\cos^{n-2}(x) \bigl[n^2 \cos^2(x) + (-n^2 + n)\bigr],\pnl$\\
$\resultleft{\cos}{3} = -\cos^{n-3}(x) \bigl[n^3 \sin^3(x) + (-3 n^2 + 2 n) \sin(x)\bigr]\pnl$\\
$\resultbox{} = \cos^{n-3}(x) \sin(x) \bigl[n^3 \cos^2(x) + (-n^3 + 3 n^2 - 2 n)\bigr],\pnl$\\
$\resultleft{\cos}{4} = \cos^{n-4}(x) \bigl[n^4 \sin^4(x) + (-6 n^3 + 8 n^2 - 4 n) \sin^2(x) + (3 n^2 - 2 n)\bigr]\pnl$\\
$\resultbox{} = \cos^{n-4}(x) \bigl[n^4 \cos^4(x) + (-2 n^4 + 6 n^3 - 8 n^2 + 4 n) \cos^2(x) + (n^4 - 6 n^3 + 11 n^2 - 6 n)\bigr],\pnl$\\
$\resultleft{\cos}{5} = -\cos^{n-5}(x) \bigl[n^5 \sin^5(x) + (-10 n^4 + 20 n^3 - 20 n^2 + 8 n) \sin^3(x)\pnl$\\
$\resultbox{}\phantom{ = -\cos^{n-5}(x)}\;+ (15 n^3 - 30 n^2 + 16 n) \sin(x)\bigr]\pnl$\\
$\resultbox{} = -\cos^{n-5}(x) \sin(x) \bigl[n^5 \cos^4(x) + (-2 n^5 + 10 n^4 - 20 n^3 + 20 n^2 - 8 n) \cos^2(x)\pnl$\\
$\resultbox{}\phantom{ = -\cos^{n-5}(x) \sin(x)}\;+ (n^5 - 10 n^4 + 35 n^3 - 50 n^2 + 24 n)\bigr],\pnl$\\
$\resultleft{\cos}{6} = \cos^{n-6}(x) \bigl[n^6 \sin^6(x) + (-15 n^5 + 40 n^4 - 60 n^3 + 48 n^2 - 16 n) \sin^4(x)\pnl$\\
$\resultbox{}\phantom{ = \cos^{n-6}(x)}\;+ (45 n^4 - 150 n^3 + 196 n^2 - 88 n) \sin^2(x) + (-15 n^3 + 30 n^2 - 16 n)\bigr]\pnl$\\
$\resultbox{} = -\cos^{n-6}(x) \bigl[n^6 \cos^6(x) + (-3 n^6 + 15 n^5 - 40 n^4 + 60 n^3 - 48 n^2 + 16 n) \cos^4(x)\pnl$\\
$\resultbox{}\phantom{ = -\cos^{n-6}(x)}\;+ (3 n^6 - 30 n^5 + 125 n^4 - 270 n^3 + 292 n^2 - 120 n) \cos^2(x)\pnl$\\
$\resultbox{}\phantom{ = -\cos^{n-6}(x)}\;+ (-n^6 + 15 n^5 - 85 n^4 + 225 n^3 - 274 n^2 + 120 n)\bigr],~\ldots\pnl$
\vspace*{-0.9ex}
\subsection*{The derivatives of the powers of hyperbolic sine}
\label{sinh}
\vspace*{-0.2ex}
$\resultleft{\sinh}{0} = \sinh^{n-0}(x) \bigl[1\bigr]\pnl$\\
$\resultbox{} = \sinh^{n-0}(x) \bigl[1\bigr] = \sinh^n(x),\pnl$\\
$\resultleft{\sinh}{1} = \sinh^{n-1}(x) \bigl[n \cosh(x)\bigr]\pnl$\\
$\resultbox{} = \sinh^{n-1}(x) \cosh(x) \bigl[n\bigr],\pnl$\\
$\resultleft{\sinh}{2} = \sinh^{n-2}(x) \bigl[n^2 \cosh^2(x) - n\bigr]\pnl$\\
$\resultbox{} = \sinh^{n-2}(x) \bigl[n^2 \sinh^2(x) + (n^2 - n)\bigr],\pnl$\\
$\resultleft{\sinh}{3} = \sinh^{n-3}(x) \bigl[n^3 \cosh^3(x) + (-3 n^2 + 2 n) \cosh(x)\bigr]\pnl$\\
$\resultbox{} = \sinh^{n-3}(x) \cosh(x) \bigl[n^3 \sinh^2(x) + (n^3 - 3 n^2 + 2 n)\bigr],\pnl$\\
$\resultleft{\sinh}{4} = \sinh^{n-4}(x) \bigl[n^4 \cosh^4(x) + (-6 n^3 + 8 n^2 - 4 n) \cosh^2(x) + (3 n^2 - 2 n)\bigr]\pnl$\\
$\resultbox{} = \sinh^{n-4}(x) \bigl[n^4 \sinh^4(x) + (2 n^4 - 6 n^3 + 8 n^2 - 4 n) \sinh^2(x) + (n^4 - 6 n^3 + 11 n^2 - 6 n)\bigr],\pnl$\\
$\resultleft{\sinh}{5} = \sinh^{n-5}(x) \bigl[n^5 \cosh^5(x) + (-10 n^4 + 20 n^3 - 20 n^2 + 8 n) \cosh^3(x)\pnl$\\
$\resultbox{}\phantom{ = \sinh^{n-5}(x)}\;+ (15 n^3 - 30 n^2 + 16 n) \cosh(x)\bigr]\pnl$\\
$\resultbox{} = \sinh^{n-5}(x) \cosh(x) \bigl[n^5 \sinh^4(x) + (2 n^5 - 10 n^4 + 20 n^3 - 20 n^2 + 8 n) \sinh^2(x)\pnl$\\
$\resultbox{}\phantom{ = \sinh^{n-5}(x) \cosh(x)}\;+ (n^5 - 10 n^4 + 35 n^3 - 50 n^2 + 24 n)\bigr],\pnl$\\
$\resultleft{\sinh}{6} = \sinh^{n-6}(x) \bigl[n^6 \cosh^6(x) + (-15 n^5 + 40 n^4 - 60 n^3 + 48 n^2 - 16 n) \cosh^4(x)\pnl$\\
$\resultbox{}\phantom{ = \sinh^{n-6}(x)}\;+ (45 n^4 - 150 n^3 + 196 n^2 - 88 n) \cosh^2(x) + (-15 n^3 + 30 n^2 - 16 n)\bigr]\pnl$\\
$\resultbox{} = \sinh^{n-6}(x) \bigl[n^6 \sinh^6(x) + (3 n^6 - 15 n^5 + 40 n^4 - 60 n^3 + 48 n^2 - 16 n) \sinh^4(x)\pnl$\\
$\resultbox{}\phantom{ = \sinh^{n-6}(x)}\;+ (3 n^6 - 30 n^5 + 125 n^4 - 270 n^3 + 292 n^2 - 120 n) \sinh^2(x)\pnl$\\
$\resultbox{}\phantom{ = \sinh^{n-6}(x)}\;+ (n^6 - 15 n^5 + 85 n^4 - 225 n^3 + 274 n^2 - 120 n)\bigr],~\ldots\pnl$
\vspace*{-0.9ex}
\subsection*{The derivatives of the powers of hyperbolic cosine}
\label{cosh}
\vspace*{-0.2ex}
$\resultleft{\cosh}{0} = \cosh^{n-0}(x) \bigl[1\bigr]\pnl$\\
$\resultbox{} = \cosh^{n-0}(x) \bigl[1\bigr] = \cosh^n(x),\pnl$\\
$\resultleft{\cosh}{1} = \cosh^{n-1}(x) \bigl[n \sinh(x)\bigr]\pnl$\\
$\resultbox{} = \cosh^{n-1}(x) \sinh(x) \bigl[n\bigr],\pnl$\\
$\resultleft{\cosh}{2} = \cosh^{n-2}(x) \bigl[n^2 \sinh^2(x) + n\bigr]\pnl$\\
$\resultbox{} = \cosh^{n-2}(x) \bigl[n^2 \cosh^2(x) + (-n^2 + n)\bigr],\pnl$\\
$\resultleft{\cosh}{3} = \cosh^{n-3}(x) \bigl[n^3 \sinh^3(x) + (3 n^2 - 2 n) \sinh(x)\bigr]\pnl$\\
$\resultbox{} = \cosh^{n-3}(x) \sinh(x) \bigl[n^3 \cosh^2(x) + (-n^3 + 3 n^2 - 2 n)\bigr],\pnl$\\
$\resultleft{\cosh}{4} = \cosh^{n-4}(x) \bigl[n^4 \sinh^4(x) + (6 n^3 - 8 n^2 + 4 n) \sinh^2(x) + (3 n^2 - 2 n)\bigr]\pnl$\\
$\resultbox{} = \cosh^{n-4}(x) \bigl[n^4 \cosh^4(x) + (-2 n^4 + 6 n^3 - 8 n^2 + 4 n) \cosh^2(x)\pnl$\\
$\resultbox{}\phantom{ = \cosh^{n-4}(x)}\;+ (n^4 - 6 n^3 + 11 n^2 - 6 n)\bigr],\pnl$\\
$\resultleft{\cosh}{5} = \cosh^{n-5}(x) \bigl[n^5 \sinh^5(x) + (10 n^4 - 20 n^3 + 20 n^2 - 8 n) \sinh^3(x)\pnl$\\
$\resultbox{}\phantom{ = \cosh^{n-5}(x)}\;+ (15 n^3 - 30 n^2 + 16 n) \sinh(x)\bigr]\pnl$\\
$\resultbox{} = \cosh^{n-5}(x) \sinh(x) \bigl[n^5 \cosh^4(x) + (-2 n^5 + 10 n^4 - 20 n^3 + 20 n^2 - 8 n) \cosh^2(x)\pnl$\\
$\resultbox{}\phantom{ = \cosh^{n-5}(x) \sinh(x)}\;+ (n^5 - 10 n^4 + 35 n^3 - 50 n^2 + 24 n)\bigr],\pnl$\\
$\resultleft{\cosh}{6} = \cosh^{n-6}(x) \bigl[n^6 \sinh^6(x) + (15 n^5 - 40 n^4 + 60 n^3 - 48 n^2 + 16 n) \sinh^4(x)\pnl$\\
$\resultbox{}\phantom{ = \cosh^{n-6}(x)}\;+ (45 n^4 - 150 n^3 + 196 n^2 - 88 n) \sinh^2(x) + (15 n^3 - 30 n^2 + 16 n)\bigr]\pnl$\\
$\resultbox{} =  \cosh^{n-6}(x) \bigl[n^6 \cosh^6(x) + (-3 n^6 + 15 n^5 - 40 n^4 + 60 n^3 - 48 n^2 + 16 n) \cosh^4(x)\pnl$\\
$\resultbox{}\phantom{ = \cosh^{n-6}(x)}\;+ (3 n^6 - 30 n^5 + 125 n^4 - 270 n^3 + 292 n^2 - 120 n) \cosh^2(x)\pnl$\\
$\resultbox{}\phantom{ = \cosh^{n-6}(x)}\;+ (-n^6 + 15 n^5 - 85 n^4 + 225 n^3 - 274 n^2 + 120 n)\bigr],~\ldots\pnl$
\section*{Source code}
\label{source}
{\sloppy Source code for a .\textsc{net} implementation, to list and calculate and evaluate this work's polynomials and expressions, written in \Fsharp\kern -0.12em,
can be found at \href{http://fssnip.net/7Xf}{\texttt{http://fssnip.net/7Xf}}\,.}
\bibliographystyle{plainnat}
\bibliography{Derivatives}
\end{document}